\documentclass[11pt]{article}

\usepackage[utf8]{inputenc}  
\usepackage[T1]{fontenc}
\usepackage{amssymb}
\usepackage{graphicx}
\usepackage{float}
\usepackage{bbold}
\usepackage{dsfont}
\usepackage{epsfig}
\usepackage{url}
\usepackage{xcolor}
\usepackage{lscape}
\definecolor{grismoi}{cmyk}{0, 0, 0, 0.65}

\begin{document}

\title{A proportional hazard model for storm occurrence risk}
\date{}
\maketitle

Malika Chassan $^{1}$ \& Jean-Marc Azaïs $^{2}$ \& Guillaume Buscarlet $^{3}$ \& Norbert Suard $^{4}$\\
\bigskip

{\it
$^{1}$ Institut de Mathématiques de Toulouse - Université Paul Sabatier, 118 Route de
Narbonne, 31062 Toulouse, Cedex 9, France. \\
E-mail : malika.chassan@math.univ-toulouse.fr
 
$^{2}$  Institut de Mathématiques de Toulouse \\
E-mail : jean-marc.azais@math.univ-toulouse.fr

$^{3}$ Thales Alenia Space, 26 Avenue Jean François Champollion, 31100 Toulouse, France

$^{4}$ CNES, 18 Avenue Édouard Belin, 31400 Toulouse, France

}

\section*{Abstract}
The aim of this paper is to give a precise estimation of the extreme magnetic storms frequency per time unit (year) throughout a solar cycle.
An innovative approach based on a proportional hazard model is developed. Based on the Cox model, this method includes non-stationarity and covariate influence. The model assumes that the number of storms during a cycle is a non-homogeneous Poisson process. The intensity of this process can be expressed as the product of a baseline risk and a risk factor. In the Cox model, the baseline risk is a nuisance parameter.  In our model, it is a parameter of interest that will be estimated. The risk factor depends on a covariate, the solar activity index. As in Extreme Value Theory (EVT) and especially in Peaks Over Threshold (POT) modeling, all the high level events are used to make estimations and the results are extrapolated to the extreme level events.
This study highlights a strong correlation between the occurrence intensity of magnetic storms and their position on the solar cycle. The model can be used to forecast occurrence intensity for the current 24th solar cycle.

\section{INTRODUCTION}
Sun activity has a direct influence on Earth ionosphere and magnetic field. Earth is protected by its magnetosphere and most of the time, solar radiations have a weak impact. However, sometimes Earth is attained by a more important solar flare or a large coronal mass ejection. These phenomena cause severe disturbances of ionosphere properties and are analysed as ionospheric or geomagnetic storms.

Strongest storms can affect the navigation applications as the SBAS (Satellite-based augmentation systems) EGNOS (European Geostationary Navigation Overlay Service), the European counterpart of the WAAS (Wide Area Augmentation System), mainly in terms of integrity and continuity. For example, the period from October 19th to November 7th 2003 was particularly perturbed and the WAAS was disabled during about 30 hours \cite{report2008}. This event, called the Halloween event, is composed of several major magnetic storms. Two of them (on 29-30 October) are among the strongest storms since 1932 and they are the first and the second largest storms of the solar cycle 23. During these two storms, the vertical position error limit for the Vertical NAVigation (VNAV) function exceeded the upper limit defined by the Federal Aviation Administration (50 meters) during more than 25 hours, making the WAAS unusable for aircraft precision approaches \cite{NOAAmemo}.

The need for investigation of extreme storm risk is demonstrated with these highly critical examples. The aim of this paper is to describe an innovative approach for the estimation of the probability of occurrence of severe ionospheric storms (per time unit (year)).

Many geomagnetic indices can be used to describe the ionosphere magnetic activity. Rifa \cite{Rifa}, shows a high correlation rate between all these indices. Thus, only one index is used to make the majority of analyses, the three-hour $ap$ index (for "planetary amplitude"). This index is obtained using measurements of magnetic field variations from 13 geomagnetic observatories over the world. Thus, the main advantage of this index is its global character. An extended description of this index and a review of the other geomagnetic indices are presented in Section \ref{section : index}.

\begin{figure}
\includegraphics[width=\textwidth]{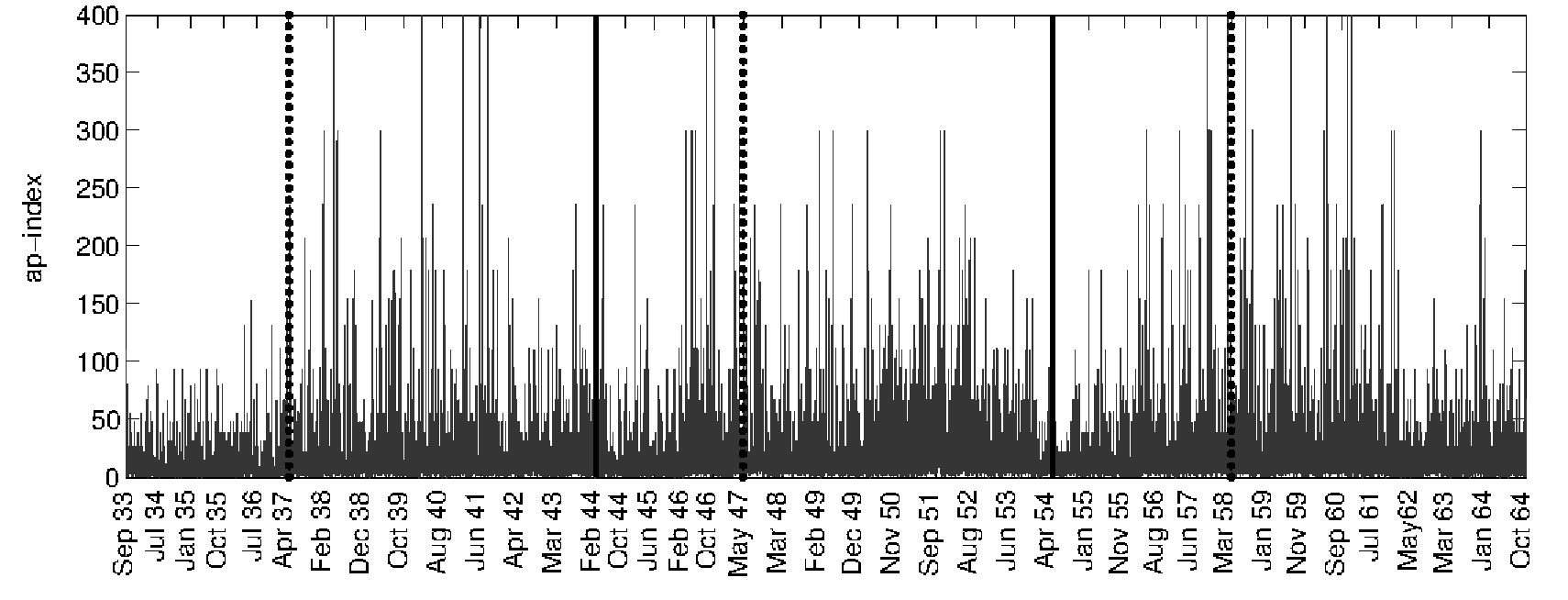} 
\caption{$ap$ index during the first three solar cycles (from the 17th to the 19th). Vertical plain lines represent the change of cycle and dotted lines correspond to the peak of each cycle.}
\label{fig : cycles123} 
\end{figure}

The data used are retrieved from the National Geophysical Data Center of the NOAA (National Oceanic and Atmospheric Administration) \cite{noaaKpap}. They consist of 80 years of registration, including 7 complete solar cycles, from the 17th (on the general list), which starts on September 1933, to the 23th, which ends on December 2008. The $ap$ index for the first three cycles of the data set is plotted in Figure \ref{fig : cycles123} on page \pageref{fig : cycles123}.

Intensive storms being scarce, classical statistical methods for probability estimation, such as empirical frequency estimation, are not accurate enough. In many domains, the probability of scarce extreme events can be estimated by extrapolation, using the Extreme Value Theory. Nevertheless, a direct application of EVT is not achievable here.

Indeed, a first issue to investigate is the non-stationarity of the studied phenomenon. Since 1843, and the publication of the short article by Schwabe \cite{Schwabe}, it is known that Sun activity (characterized by the number of solar sunspots) is characterized by 11-year cycles. Corresponding cycles are observable in the $ap$ index data set (Figure \ref{fig : cycles123}). Hüsler \cite{Husler86} extends the results of EVT for stationary processes to some non-stationary ones but no general theory exists for non-stationary processes (\cite{Coles2001}, Chapter 6).

However, various non-stationary extreme models are proposed, as the pragmatic approach of Coles (again Chapter 6 in \cite{Coles2001}) or the work of Jonathan and Ewans \cite{Jonathan2011}. In this paper, authors model the seasonality of extreme waves in the Gulf of Mexico. Occurrence rate and intensity of storm peak events vary with season (the Generalized Pareto parameters are expressed as a function of seasonal degree using a Fourier form). Our model could have adopted these approaches but the following arguments show that it was not relevant.


The main problem to the classical application of EVT comes from the type of data used for this study. The $ap$ index lies in the set $\{$0, 2, 3, 4, 5, 6, 7, 9, 12, 15, 18, 22, 27, 32, 39, 48, 56, 67, 80, 94, 111, 132, 154, 179, 207, 236, 300, 400$\}$. Thus, it is a discrete variable, taking a finite number of values. On one hand, the discrete case is explored in \cite{Anderson1970} where Anderson provides an extension of classical convergence results to some discrete variables. However, in practice, working with a discrete support distribution can be troublesome in some applications. As example, Cooley, Nychka and Naveau \cite{CooleyNaveau2007} analyse hydrometric data but they found that the low precision of the measurements makes the data almost discrete. Hence, a sawtooth behavior of the parameter estimators  was found with increasing threshold, which hampered to clearly select a value for the threshold.
On the other hand, the fact that the $ap$ index takes a finite number of values seems to be incompatible with the application of EVT. Indeed, an assumption used by Hüsler \cite{Husler86} in the study of exceedances behavior is violated when considering a variable with a finite support. Similarly, the result obtained by Anderson \cite{Anderson1970} assumes that the distribution support is unbounded. Hence, to the best of our knowledge, the application of EVT to the $ap$ index data set is not achievable.


Thus, an innovative approach is developed to estimate the occurrence probability of severe magnetic storms. This model, based on proportional hazard modeling, is also based on EVT. As in the classical Peaks Over Threshold method, the estimation results are obtained using high and less rare events and extrapolated to extreme events.

The main contribution of this paper is the description of a new proportional hazard model, which is detailed in Section \ref{section : model}, along with the Cox method. Definitions of the handled objects and data pretreatments can be found in Section \ref{section : data}. The descriptions of parameter estimators are gathered in Section \ref{section : estimation}. Section \ref{section : results} is dedicated to data applications and method extensions.

\section{MODEL DESCRIPTION}
\label{section : model}
The approach developed in this paper is based on the Cox proportional hazard model. Firstly, the general notion of proportional hazard model is presented and details on the particular Cox model are given. Secondly, the counting process formulation of these models is introduced. Finally, our model is described with all the modifications from the initial Cox model.

\subsection{Cox model and Proportional hazard}
\label{section : survie}
First introduced in epidemiology \cite{Cox72}, the Cox model is a survival model. Survival analysis consists in considering an individual (or a sample of individuals) from a starting time ($t_0=0$ without loss of generality) and observing it to detect the occurrence of an event of interest that might happen once. 

The occurrence time can be represented by a random variable, $T$, which also corresponds to survival time (since $t_0=0$). The variable $T$ is used to define the instantaneous risk (or hazard rate) for an individual. The hazard rate at time $t$, $\lambda(t)$, is the probability that an event occurs in the infinitesimal interval $[t, t+dt]$ for an individual who has not experienced the event until $t$:
%
\begin{equation}
\lambda(t) = \lim_{dt \to 0} \frac{1}{dt} \mathds P( t\leq T <t+dt | T  \geq t),
\label{eq : def_lambda}
\end{equation} 
where $\mathds P(A)$ is the probability of the event $A$. 

The Cox model expresses the instantaneous risk with respect to time $t$ and the covariates $(X_1,...,X_p)$ as:
\begin{equation}
\lambda(t,X_1, ..., X_p) = \lambda_0(t)\exp (\sum_{i=1}^p \beta_i X_i),
\label{eq : cox}
\end{equation}
where $\lambda_0(t)$ is the baseline risk.

A consequence of the model (\ref{eq : cox}) is that the instantaneous risk is proportional to a baseline risk and that the coefficient is the exponential of a linear combination of the covariates.
%
%
%
%
Thus, considering two individuals with covariates $(X^{(1)}_1, ..., X^{(1)}_p)$ and $(X^{(2)}_1, ..., X^{(2)}_p)$ the risk factor for the first patient, compared to the second, only depends on covariate effects:

$$\frac{\lambda(t,X^{(1)}_1, ..., X^{(1)}_p)}{ \lambda(t,X^{(2)}_1, ..., X^{(2)}_p) } = \exp \left(\sum_{i=1}^p \beta_i (X^{(1)}_i-X^{(2)}_i)\right).$$

In survival analysis, data used are often censored. 
In this paper, data are uncensored since only complete solar cycles are considered. Hence, the censoring question will not be developed further. For more details about survival analysis and proportional hazard models, see Aalen, Borgan and Gjessing \cite{Aalen2008}. 

\subsection{Counting process}
\label{section : counting}
Now, consider a recurrent event that can happen several times for the same individual and the counting process $N(t)$ that counts the number of events up to $t$. In analogy with Formula (\ref{eq : def_lambda}), we define a proportional intensity model by:
$$
\mathds P( N(t+dt)-N(t) = 1) = \lambda_0(t)\exp (\sum_{i=1}^p \beta_i X_i)dt.
$$
Andersen and Gill, in \cite{AndersenGill}, define a counting process formulation of the Cox model, including the possibility for the event to be recurrent.

\subsection{Our model}
In this paper, a proportional hazard model is applied. The event of interest is the occurrence of a storm and a solar cycle corresponds to an individual (indexed by the subscript $j$). Only one covariate is considered, the solar activity index of the cycle $j$, $X_j$, and the parameter $\beta$ models its influence. The model constructed has undergone meaningful modifications from the Cox model:

- observations are not censored;

- the number of events (storm occurrences) up to time $t$, $N(t)$, is considered. Hence the counting process is supposed to be an inhomogeneous Poisson process with intensity depending on time and covariates;

- the variable $D_j$ (length of cycle $j$, in years) is included as factor. Thus, the measurement unit is the number of events per time unit and not per cycle;

-  $\lambda_0(t)$ is not considered as a nuisance parameter but as a parameter to estimate;

- as in POT modeling in EVT, estimation is made using all the high level storms and an extrapolation to the storms of extreme level is applied using the parameter $P_{400}$, the probability that a high level storm grows into a storm of level 400 \footnote{400 is the value of the $ap$ index characterizing extreme storms. More details about high and extreme levels are given in Sections \ref{section : storm def} and \ref{section : data sum}}. Utilization of this parameter assumes that the level reached by a high storm does not depend on the instant of appearance. A chi-square independence test showed that this assumption is acceptable (Appendix \ref{annexe: chi2}).

Thus, in our model, the number of observed storms (of high level) during cycle $j$ up to time $t$, called $N_j(t)$, is supposed to be a non-homogeneous Poisson process of intensity $\lambda_j(t)$ such as:
$$ \lambda_j(t) = \lambda_0(t) D_j \exp ( \beta X_j),$$
i.e.
$$  N_j([a,b]) \sim \mathcal P \left(  \int_a^b \lambda_j (t)dt  \right),$$
with $\mathcal P(\xi)$ denoting the Poisson distribution of intensity $\xi$.

The baseline intensity $\lambda_0(t)$ takes into account the fact that storms occur more likely during the second half of the cycle. It is the main parameter.

Throughout this paper $\lambda_0(t)$ represents the baseline intensity of high level storms, $\lambda_j(t)$ represents the instantaneous intensity of high level storms, which takes into account the risk factor $D_j \exp ( \beta X_j)$ (for cycle $j$). The extrapolated instantaneous intensity for extreme level storms (level 400) is  $\lambda_{400, j}(t)$, i.e. $\lambda_{400, j}(t) = \lambda_j(t) \times P_{400}$. 
Unless the indices j are specified, the parameters refer to a theoretical cycle with a mean solar activity (146.7 in the case of the $ap$ data set).

\section{DATA AND PRETREATMENTS}
\label{section : data}
Prior to model parameter estimation, in this Section, the used data and their pretreatments are illustrated. The choice of the $ap$ index is discussed and definitions of a storm as well as high and extreme levels are given.


\subsection{Geomagnetic indices}
\label{section : index}

The most commonly used geomagnetic indices are listed below, each with its advantages and disadvantages:

-  the $K$ index, introduced by Bartels in 1938, is given with a sampling period of three hours. Each measurement corresponds to the maximal variation of the horizontal component of Earth magnetic field (sum of maximum positive and negative deviations, in nano-Tesla) during the past three-hour interval, compared to a reference quiet day. Each value of the $K$ index is a digit between 0 and 9. As its quasi logarithmic scale  limits the use of the $K$ index, especially for calculation of averages, additional indices were created.

- The $a$ index (for "amplitude") is a linear transformation of $K$. The $A$ index is the average of the 8 $a$-values per day. The $K$, $a$ and $A$ indices are computed by various observatories around the world and measurements have to be merged to describe the global behavior of ionosphere.

- Local variations are smoothed by considering weighted averages over the world, represented by the $ap$ and $Ap$ indices, where $p$ stands for "planetary" and the capital letter $A$ refers to the daily average. These two indices are computed using the $K$ index values from 13 observatories distributed over Earth. Weights are fixed according to observatories latitudes. A planetary version of $K$ also exists ($Kp$) and Table \ref{tab : relationKap} gives the relation between $ap$ and $Kp$. The $Kp$ and $ap$ indices have been available since 1932.


- The $aa$ index (for "antipodal amplitude") is computed using the $K$ index from only two nearly antipodal geomagnetic stations located in England and Australia. Thus, this index can not be considered as global as the $ap$ index. However, it has the advantage of being available since 1868.

\begin{figure}[H]
\center
\includegraphics[width=8cm,height=5.5cm]{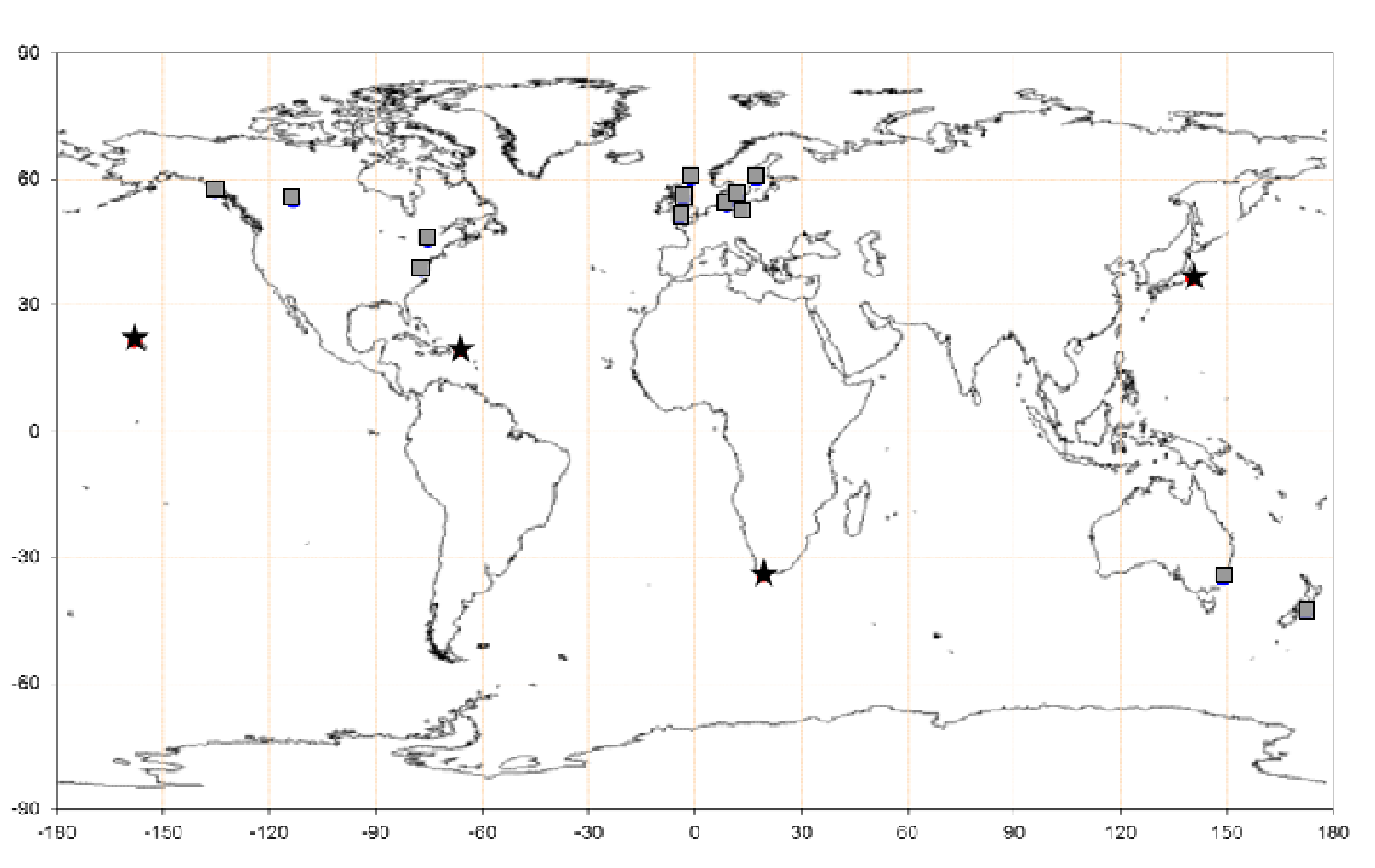} 
\caption{Positions of observatories for the $Dst$ ($\bigstar$) and the $Kp$/$ap$ indices  (\textcolor{gray}{$\blacksquare$}).}
\label{fig : carte_obs} 
\end{figure}

- The $Dst$ (Disturbance storm time) index is restricted to the equatorial magnetic perturbation (see Figure \ref{fig : carte_obs}). 57 years of data are available (against 80 for the $ap$ index). 

- The raw geomagnetic data are available for many geomagnetic observatories. Oldest data date back to 1883 for hourly values and to 1969 for 1 minute values. They consist of measurements of magnetic field variations by means of magnetometers. One disadvantage of these data is the presence of recording gaps (with gap lengths varying from one month to several years depending on the observatory). However, the principal drawback of these data is the level of required pre-treatments.
See \cite{noaaKpap} for more details on geomagnetic indices.

The $ap$ data set has the advantage of being larger than the $Dst$ one and it does not exhibit recording gaps, contrary to raw geomagnetic data. Moreover, the chosen $ap$ index is more global than the $aa$ or $Dst$ indices, as can be seen in Figure \ref{fig : carte_obs}.

\subsection{Declustering and storm definition}
\label{section : storm def}
Ionospheric perturbations are classified in a standardized way using the $ap$ index, according to Table \ref{tab : relationKap} \cite{Norbert}:
\begin{table}[H]
\centering
\begin{tabular}{|l|c|c|}
\hline
\textbf{Ionosphere Condition}  & \textbf{ $Kp$ index} &  \textbf{$ap$ index} \\ \hline  \hline
Quiet &   0-1 & <7   \\ \hline
Unsettled  & 2 &  7 to <15     \\ \hline
Active  & 3 &  15 to <27       \\ \hline
Minor storm & 4 &  27 to <48   \\ \hline
Major storm & 5 &    48 to <80 \\ \hline
Severe storm & 6 &  80 to <132 \\ \hline
Large severe & 7 &  132 to <236 \\ \hline
\hspace{0.4cm}Extreme & 8 &  236 to <400\\ \hline
\hspace{0.4cm}Extreme & 9 &   400 \\ \hline
\end{tabular}\\
\caption{Relation between $Kp$, $ap$ and ionosphere activity.}
\label{tab : relationKap}
\end{table}

In this work, a storm is defined by the event with the highest intensity even though a storm normally presents different periods  of high intensity separated by less active ones. In order to detect the event of highest intensity, a declustering process (referred to as "runs declustering" in \cite{Coles2001}, Chapter 5.3) is applied. Two parameters have to be set:

- a low level, the threshold above which a storm is considered to begin (typically 111, 132 or 154);

- the run length \textit{r} which is the minimal number of observations below the low level between two events for them to be considered as independent. 

Thus, two exceedances of the low level separated by less than \textit{r} measurements are considered to belong to the same cluster (same storm). For each cluster, the storm level is defined as the maximal level reached in the cluster. The first occurrence of this maximum is also recorded, representing the storm date. For a cluster, the storm length is the number of observations between the first up-crossing and the last down-crossing of the low level.

Durations of magnetic storms are highly variable, from 3 or 6 hours for a very high or extreme storm (level 300 or 400) to 90 hours for a low level storm. Due to the declustering in this work, only one-time events are considered (only the first maximum occurrence time is recorded). This is not incoherent since our analyses focus on strong storms, which are brief compared to longer lasting but weaker storms. Nevertheless, the duration of strong storms should be taken into account for the probability of occurrence definition (see Section \ref{section: precise_proba}).

\subsection{Data summary}
\label{section : data sum}
The $ap$ data set contains 7 complete solar cycles and after declustering only 23 magnetic storms of level 400 are detected. The number of events is too low to estimate their frequency as a function of covariate. For storms of level 300, 44 events are detected and this is still insufficient for frequency estimation. 

Consequently, as in Extreme Value Theory, less rare events are used to estimate influence of each covariate and an extrapolation of the results to extreme level is performed. In this paper, less rare events are the high level storms (but not necessary extreme). Thus, all the storms of level greater or equal to the low level parameter defined in the declustering process are used. For example, if the low level is 111, a "high level storm" is any storm of level 111, 132, 154, 179, 207, 236, 300 or 400. The "extreme level" is only 400. 
Mean probability of occurrence for each high level is given in Table \ref{tab: freq}.

\begin{table}[H]
\centering

\begin{tabular}{|l|c|c|c|c|c|c|c|c|}
\hline
Level  & Number of storms & Freq.($\times 10^{4}$) & Freq.(year$^{-1}$)\\ \hline
\hline
111 & 182 & 7.99 & 2.33\\ \hline
132 & 158 & 6.93 & 2.02\\ \hline
154 & 103 & 4.52 & 1.32\\ \hline
179 & 84 & 3.69  & 1.08\\ \hline
207 & 51 & 2.24  & 0.65\\ \hline
236 & 57 & 2.50  & 0.73\\ \hline
300 & 44 & 1.93  & 0.56\\ \hline
400 & 23 & 1.01  & 0.29\\ \hline
\end{tabular}
\caption{Number of occurrences and frequency of storms by level.}
\label{tab: freq}
\end{table}

Besides the three-hour $ap$ index, a covariate representing the solar activity of a cycle is available and is denoted as $X$ (the unity is the number of sunspots). The solar activity of a cycle is the maximum of monthly Smoothed Sunspot Number (monthly SSN). For an easier interpretation of results, this covariate is centered. See \cite{nasaSSN} for more details on sunspot number. For the seven cycles of the data set, the solar activity is $X = ( 119.2,\  151.8,\  201.3,\  110.6,\  164.5,\  158.5,\  120.8 )$ with a mean of 146.7.

For each cycle, beginning and end dates are also available, as well as peak date. Peak date corresponds to the date of maximal solar activity within a cycle (characterized by the highest sunspot number). Lengths of cycles are gathered in the vector $D$, with $D_j$ being the length of the cycle $j$. All the covariates for the $ap$ data set are summarized in Table \ref{tab: covar}. In addition, beginning date of the current solar cycle (24th) can be found in the last line of Table \ref{tab: covar}. Values of end and peak dates as well as length and solar activity index of this 24th cycle are predictions (in \textcolor{grismoi}{gray}), since it is not ended. End of the 24th cycle is estimated around \textcolor{grismoi}{December 2019}. The NOAA values are used \cite{noaaPrev24}. 

\begin{table}[H]
\centering
\begin{tabular}{|c|c|c|c|c|c|}
\hline
Cycle & Beginning & Peak  & $X$ & $D$ & Extreme \\ 
& &  & & & storms \\ \hline
\hline
17 & Sep 1933 & Apr 1937 &  119.2 & 10.4 & 5\\ \hline
18 & Feb 1944 & May 1947 &  151.8 & 10.2 & 2\\ \hline
19 & Apr 1954 & Mar 1958 &  201.3 & 10.5 & 7\\ \hline
20 & Oct 1964 & Nov 1968 &  110.6 & 11.7 & 3\\ \hline
21 & Jun 1976 & Dec 1979 &  164.5 & 10.3 & 2\\ \hline
22 & Sep 1986 & Jul 1989 &  158.5 & 9.7 & 1\\ \hline
23 & May 1996 & Mar 2000 &  120.8 & 12.6 & 3\\ \hline
24 & Dec 2008 & \textcolor{grismoi}{Nov 2013} & \textcolor{grismoi}{87.9} &\textcolor{grismoi}{11} & \textcolor{grismoi}{-}\\ \hline
\end{tabular}\\
\caption{Beginning and peak dates, solar activity index $X$ (in sunspots), length $D$ (in years) and number of extreme storms for each cycle from the 17th to the 24th (predicted values are written in \textcolor{grismoi}{gray}).}
\label{tab: covar}
\end{table}

\subsection{Time Warping}

\begin{figure}[H]
\center
\includegraphics[width=8cm,height=6cm]{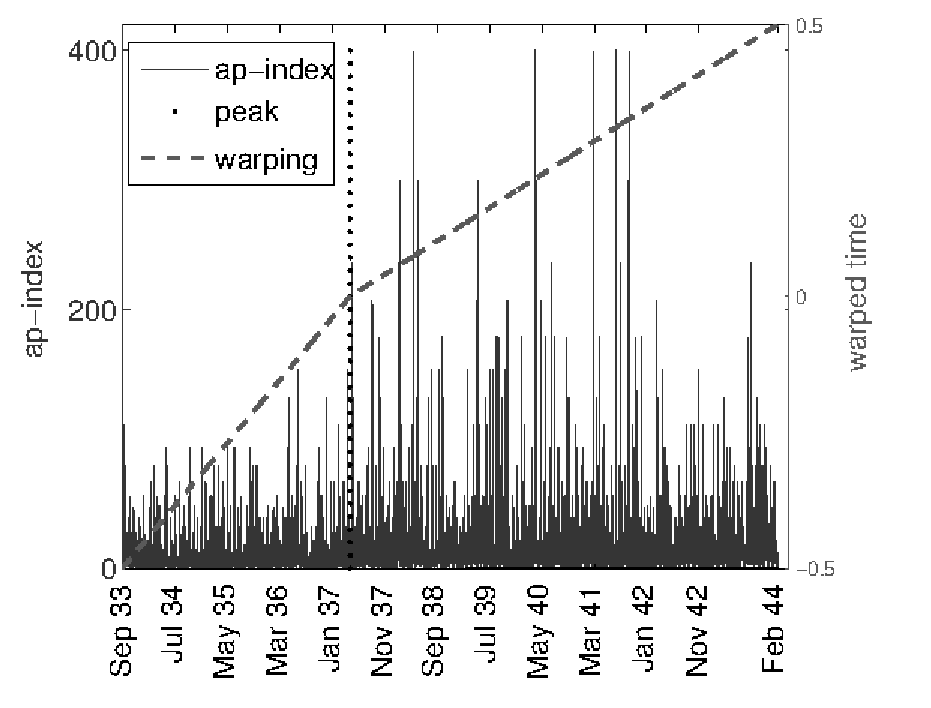} 
\caption{The $ap$ index during the first cycle (from September 1933 to February 1944). The dotted vertical line represents the peak and the dash-dotted line the warped time.} 
\label{Cycle2Warping} 
\end{figure}

Durations of the 7 complete solar cycles range from 9.7 to 12.6 years. Thus, in order to analyse all the 7 cycles together, a time warping is applied to each cycle: the position of a storm within a cycle is represented by a number between $-0.5$ and $0.5$ where  $-0.5$ is the beginning of the cycle, 0.5 its end and 0 its middle (peak). In Figure \ref{Cycle2Warping}, the dash-dotted line represents the warped time for the first complete solar cycle.

Notice that the addition of the cycle's length as factor in $\lambda_j(t) = \lambda_0(t) D_j \exp ( \beta X_j),$ assures that the measurement unit is the number of storms per year and not per cycle. Thus, the time standardization may not be biasing the intensity estimation.
In the rest of this paper, unless otherwise indicated, warped time is considered.

\subsection{Precision on probability of occurrence}
\label{section: precise_proba}

Due to the application of a declustering process, a storm is now defined by three values: the maximal level, the first time when this maximum is attained and the length of the cluster. 
Using this modeling, we can estimate the probability: 
$$  P_1(t) = \mathds P(\textrm{a storm of level 400 \textbf{begins} at time t}),$$
and we want to know the probability:
$$  P_2(t) = \mathds P(\textrm{a storm of level 400 \textbf{ is ongoing} at time t}).$$

In the whole data set, the level 400 is reached 29 times, but only 23 storms of level 400 are counted after declustering. Among these 23 storms, 17 reach the level 400 only one time and 6 remain at this level two consecutive measurement times. 

So, assuming that duration of an extreme storm is never exceeding 6 hours, and denoting $p = \mathds P(\textrm{storm stays at the level 400 two times})$, $P_2 (t)$ can be expressed as:
$$
P_2(t)  \simeq P_1(t) + P_1(t-1)\times  p.
$$
Using a local stationarity argument, we can declare that $P_1(t-1) \simeq P_1(t)$. Thus, an approximation of $P_2 (t)$ is:
$$
P_2(t) \simeq  P_1(t) \times (1 +  p),
$$
and using the $ap$ data set, this is estimated by:
$$P_2(t) \simeq  P_1(t)\times (1 +  6/23).$$

\section{ESTIMATION}
\label{section : estimation}
\subsection{$P_{400}$ and $\beta$}
Since $P_{400}$ is independent of the position within the cycle (Appendix \ref{annexe: chi2}), the empirical frequency is used
$$ \widehat{P_{400}} = \frac{\# \{ \textrm{storms of level } 400 \}}{\# \{  \textrm{storms of level $\geq$  low level} \}}.$$
A $95\%$ confidence interval is computed via Gaussian approximation. With $m=\# \{  \textrm{storms of level $\geq$ low level}\}$, the number of storm with a level above the low level, we have: 
$$P_{400} \in \left[  \widehat{P_{400}} \pm 1.96 \sqrt{\widehat{P_{400}}(1- \widehat{P_{400}})/m}\right].$$
For $\beta$, we use the fact that
$$N_j = N_j([-0.5,0.5]) \sim \mathcal P \left(  \left[ \int_{-1/2}^{1/2} \lambda_0 (s)ds \right]  \ D_j \exp(\beta X_j) \right).$$
As in the Cox model, the sufficiency of the statistic $N_j$ is verified and $\beta$ is estimated by its maximum likelihood estimator in a Poisson generalized linear model. A confidence interval is also computed. More details can be found in Appendix \ref{annexe}.
\subsection{Baseline intensity $\lambda_0(t)$}
\label{section : lambda 0 chap}

For the estimation of the baseline intensity, a kernel estimator is used. The kernel estimator used here is a slightly modified version of the extension of the Hjort's estimator, proposed by Nielsen and Linton in \cite{NielsenLinton}. In our model, the covariate does not depend on time. Thus, some terms of the Nielsen and Linton's estimator become constant. All constant terms are gathered in a normalization constant. Another difference is the fact that in our approach, the parameter of interest is $\lambda_0(t)$ instead of $\lambda(t)$, but the normalization constant is set in order to have a consistent estimator.

Assuming that $\beta$ is known, the estimator of $\widehat{\lambda_0 (t)}$ is:
$$
\begin{array}{ll} 
\widehat{\lambda_0 (t)} &= C \displaystyle \sum_{j=1}^J  \int_{-1/2}^{1/2} dN_j(t-s)\phi(s) \\
&= C \displaystyle \sum_{j=1}^J  \int_{-1/2}^{1/2} N_j(t-s)\phi'(s)ds,
\end{array}
$$
where $J$ is the number of individuals (cycles), $C$ a normalization constant  and  $\phi$ the kernel, verifying $\phi(\pm 1/2) =0$ (for the integration by parts) and $\int_{-1/2}^{1/2} \phi(s)ds =1$.\\

Bias and variance of this estimator are computed using step functions and by passage to limit. Let $\phi$ be a step function, 
$$
\phi(s) = \sum_{i=1}^n a_i \mathds 1 _{A_i}(s),
$$ 
where $\mathds 1 _{A_i}$ is the indicator function of the interval $A_i$, the $A_i = [ t_i, t_{i+1} ] $ form a partition of $[ -1/2,1/2 ]$ (we can assume $t_i<t_{i+1}$ without loss of generality) and $a_i$ are such that $\int_{-1/2}^{1/2} \phi(s) ds =1$.
Then, for each $t \in [-1/2,1/2]$,
$$
\begin{array}{ll}
\widehat{\lambda_0 (t)} &= \displaystyle \sum_{j=1}^J C \int_{-1/2}^{1/2} dN_j(t-s)\phi(s)\\
&= C \displaystyle{ \sum_{j=1}^J}  \bigg\{ a_ 1 N_j([t-t_2,t-t_1]) + ...\\
& \qquad \qquad \qquad + a_n N_j([t-t_{n+1},t-t_n]) \bigg\}.
\end{array}
$$
Thus, since $N_j([a,b]) \sim \mathcal P \left( Q_j \, \int_a^b  \lambda_0 (s)ds \right) \ $ with $Q_j = D_j \exp(\beta X_j)$ and since $\mathds E \mathcal P (\xi) = \mathds V \mathcal P (\xi) = \xi$, we get:
\begin{equation}
\mathds E \, \widehat{\lambda_0 (t)} =  C \displaystyle{ \sum_{j=1}^J} Q_j \int_{-1/2}^{1/2} \lambda_0 (s)\phi(t-s)ds,
\label{eq : esperance} 
\end{equation}
where $\mathds E$ is the expected value and $\mathds V$ the variance.
Similarly, for the variance:
\begin{equation}
\mathds V \, \widehat{\lambda_0 (t)} 
=C^2 \displaystyle{ \sum_{j=1}^N} Q_j \int_{-1/2}^{1/2} \lambda_0 (s)\phi^2(t-s)ds.
\label{eq : var}
\end{equation}
Equations (\ref{eq : esperance}) and (\ref{eq : var}) are established for step functions. Extension to general functions is straightforward by a monotone convergence argument.

In case of a kernel concentrated around zero we obtain:
$$
\mathds E \, \widehat{\lambda_0 (t)} \simeq   C \displaystyle{ \sum_{j=1}^J} Q_j \lambda_0 (t).
$$
Hence, the choice  $C = 1/ \sum Q_j$ is convenient and we get:
$$
\mathds V \, \widehat{\lambda_0 (t)}\simeq \frac{1}{ \sum Q_j}\lambda_0(t) \int_{-1/2}^{1/2}\phi^2(s)ds.
$$
In this paper, a Gaussian kernel is used, i.e.
$$
\phi(s) = \frac{1}{\sqrt{2\pi}h}\exp(-\frac{s^2}{2h^2}),
$$
with $h$ being the band-width parameter, determined later. Then, using the fact that 
$$
\phi^2(s) = \frac{1}{2\sqrt\pi h} \phi(\sqrt 2 s),
$$
with $\phi(\sqrt 2 s)$ being the density function of a normal distribution $\mathcal N \left( 0, (h/\sqrt 2 )^2 \right)$ and assuming sufficiently small $h$ values, it follows that:
$$
\int_{-1/2}^{1/2}\phi^2(s)ds \simeq \int_{-\infty}^{+\infty}\phi^2(s)ds =  \frac{1}{2\sqrt\pi h}.
$$
To avoid edge effects, that make the estimated intensity tending artificially to zero at $t = \pm 0.5$, a periodization of data is applied before estimation process. The band-width parameter $h$ is chosen by cross-validation with minimization of Integrated Square Error. See \cite{Bowman1984} or \cite{Hall1983} for more details.

Note that, as indicated in Section \ref{section: precise_proba}, the intensity estimated by kernel method does not correspond to the intensity we want to evaluate. Instead, it corresponds to the probability $P_1$ that a storm of level 400  begins at time $t$. Hence a correction is applied multiplying by $\widehat {\lambda_0(t)}$ by 29/23.

Thus, the approximate confidence interval for $\lambda_0(t)$ is:
$$
\lambda_0(t) \in \left[ \widehat{\lambda_0 (t)} \pm 1.96 \sqrt{\frac{1}{\sum Q_j}\frac{\widehat{\lambda_0 (t)}}{2\sqrt \pi h}}  \right].
$$
We recall that $Q_j = D_j \exp(\beta X_j)$. Thus, in practice, a plug-in estimator of $Q_j$ is computed, using $\hat \beta$.

\subsection{Simulations}

Simulations are performed in order to assess the quality of estimation method. Sample size ($J = $ number of cycles used) varies from 2 to 300. For each value of $J$, errors are computed over 700 estimation results. 
For a fixed number of cycles, a simulation is performed as follow:\\
- Firstly, $X$, $D$ and the parameter $\beta^*$ are generated according to values of the $ap$ data set. For example, for $\beta^*$, we use a normal distribution $\mathcal N (0.006, 0.002^2)$.\\
- Secondly, an intensity function $\lambda_0^*$ is created using Gaussian distributions and periodic functions.\\
- Finally, for each cycle, an inhomogeneous Poisson process of intensity $\lambda_0^*(t) D_j \exp(\beta^* X_j)$ is simulated and estimations are performed using our model.

\begin{figure}[H]
\includegraphics[width=6.9cm, height=6.5cm]{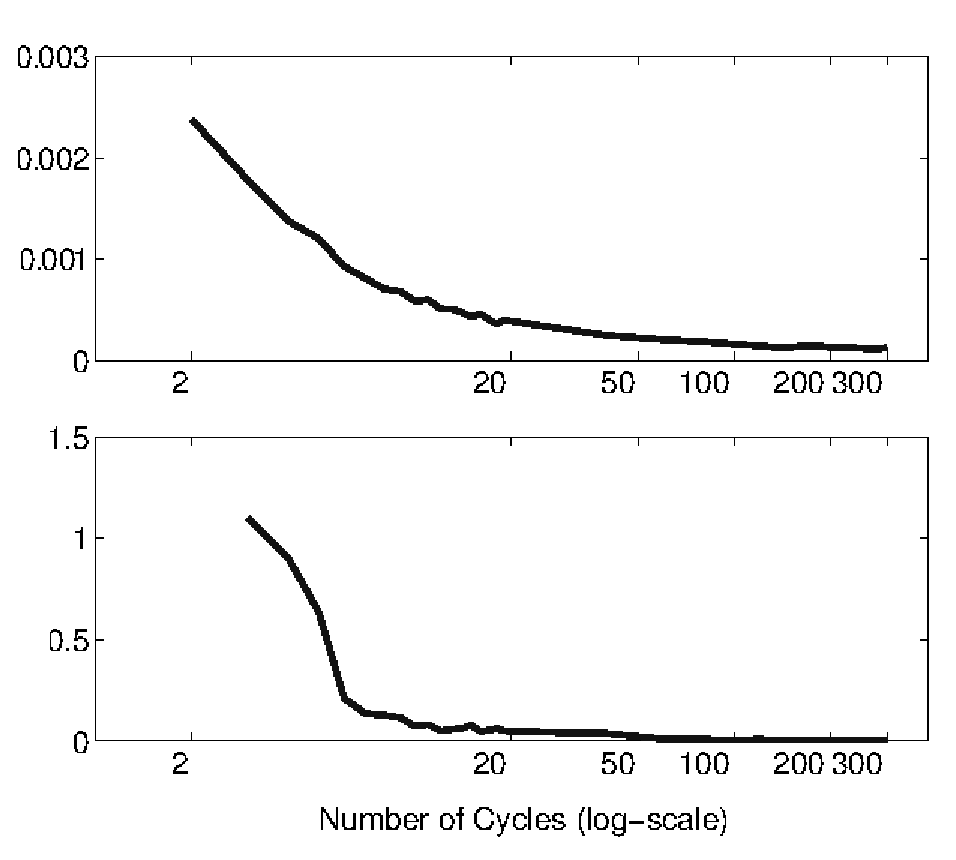}
\caption{Mean relative errors of estimation for $\lambda_0(t)$ (top) and $\beta$ (bottom) as functions of the number of cycles (in logarithmic scale). Mean values are computed over 700 error values.}
\label{fig : err_simu}
\end{figure}

Then, for each value of $J$, mean relative errors $ \mathds E \left( \frac { \hat \beta - \beta^*}{\beta^*} \right) ^2 \textrm{ and } \, \,  \mathds E \left( \frac {|| \widehat{\lambda_0} - \lambda_0^* ||_2  }{||\lambda_0^*||_2 }  \right) ^2 $ are computed and plotted in Figure \ref{fig : err_simu}.

Relative square error for $\beta$ is larger than for $\lambda_0(t)$ but it decreases more rapidly for small values of $J$. Before $J$ = 6, the relative square error is on average greater than $30 \%$. Indeed, for a very small number of cycles, the maximum likelihood estimator of $\beta$ has a large variance and estimation is not very accurate. However, for 7 cycles (as in the $ap$ data set) the mean relative square error is about $0.14$, which is acceptable.

\begin{figure}[H]
\includegraphics[width=8cm, height=6.3cm]{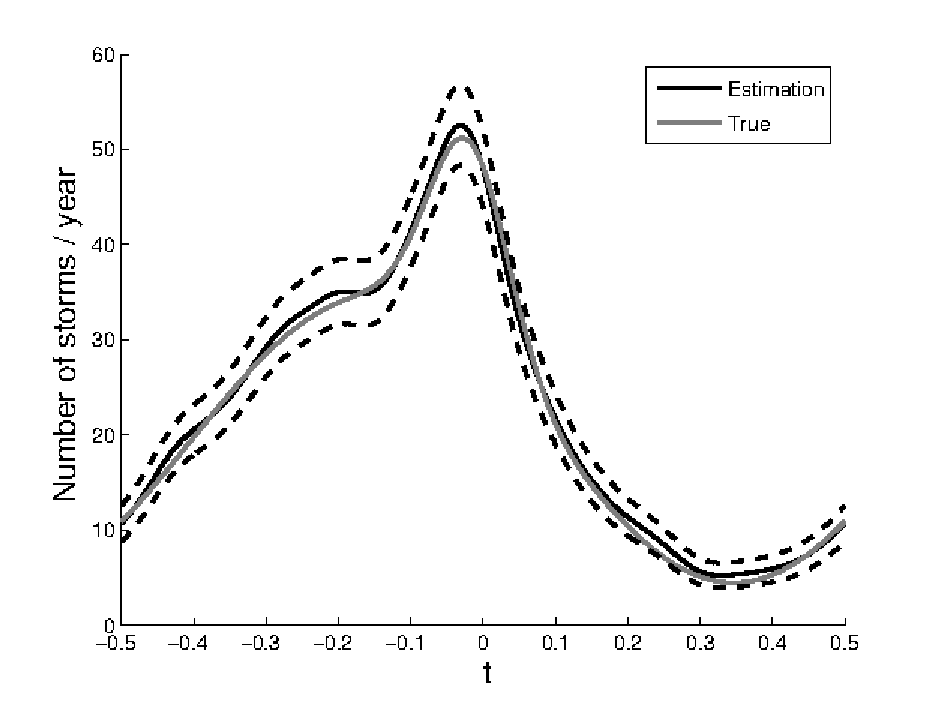}
\caption{Estimation of $\lambda_0(t)$ with a sample of 10 cycles. The curve to estimate is plotted in gray.}
\label{fig : simu}
\end{figure}

In Figure \ref{fig : simu}, an example of the baseline intensity estimated with our model and using a sample of 10 cycles is illustrated.

In accordance with what is shown in Figure \ref{fig : err_simu}, one can see that kernel estimation of baseline intensity is rather accurate. Other kernels are tested (uniform, triangular and Epanechnikov) and results are not sensitive to kernel change, for both simulated intensities and intensities estimated from real data.

\section{RESULTS} 
\label{section : results}
\subsection{Baseline intensity $\widehat {\lambda_0(t)}$}

\begin{figure}[H]
\center
\includegraphics[width=8cm,height=6.3cm]{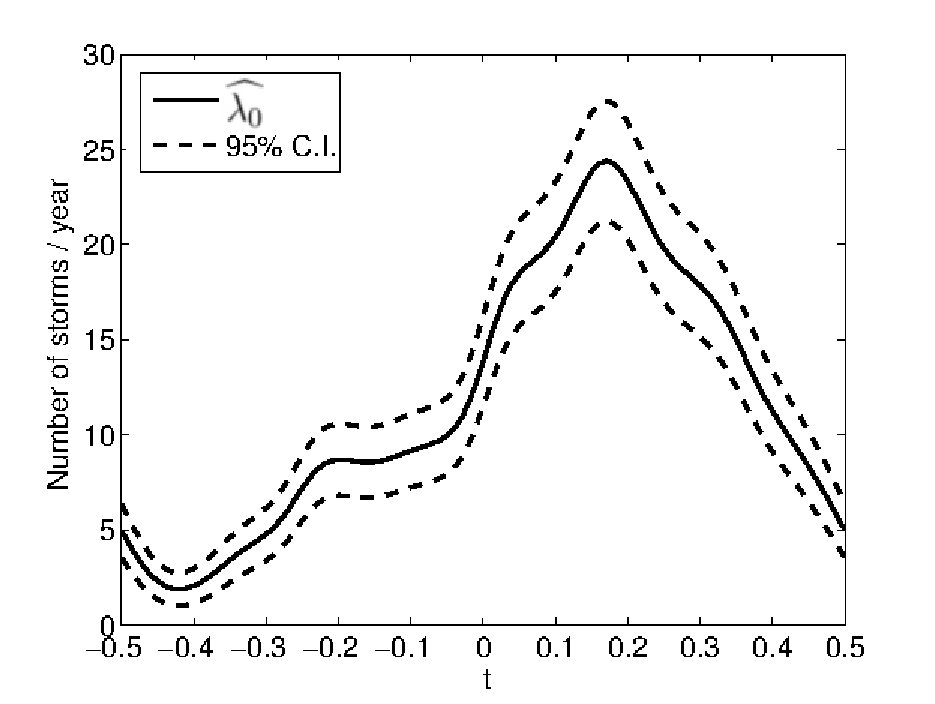}
\caption{Estimated baseline intensity (years$^{-1}$), with confidence interval, of the storms of level greater or equal to 111, for a mean solar activity of 146.7.} 
\label{fig: lambda0_111} 
\end{figure}

Figure \ref{fig: lambda0_111} depicts the estimation of $\widehat {\lambda_0(t)}$ for a low level of 111 and a run length $r = 7$  including the confidence area (i.e. the intensity for all the storms of level greater or equal to 111, for a theoretical cycle with $X=146.7$). The band-width parameter is selected by cross validation and is equal to 0.035. The baseline intensity is higher during the second half of cycle with a significant increase around $t=0$, highlighting the difference between the two halves of a solar cycle. This behavior is coherent with empirical observations of frequency of geomagnetic storms (Figure \ref{fig : cycles123}).

\subsection{$P_{400}$ and $\beta$}

For $\widehat{P_{400}}$, the results for different low levels are gathered in Table \ref{tab: p400}.

\begin{table}[H]
\center
\begin{tabular}{|l|c|c|c|}
\hline
Low level  & 111 & 132 & 154\\ \hline
\hline
$\widehat{P_{400}}$ in $\%$ & 3.14 &   4.19   &   5.93   \\ \hline
95 \% C.I. & [1.85 ; 4.43] & [2.48 ; 5.90] & [3.53 ; 8.33]  \\ \hline
\end{tabular}\\
\caption{$\widehat{P_{400}}$ (probability for a high storm to grow into a storm of level 400) and 95 \% confidence interval for each low level.}
\label{tab: p400}
\end{table}

\begin{figure}[H]
\center
\includegraphics[width=8cm,height=6.3cm]{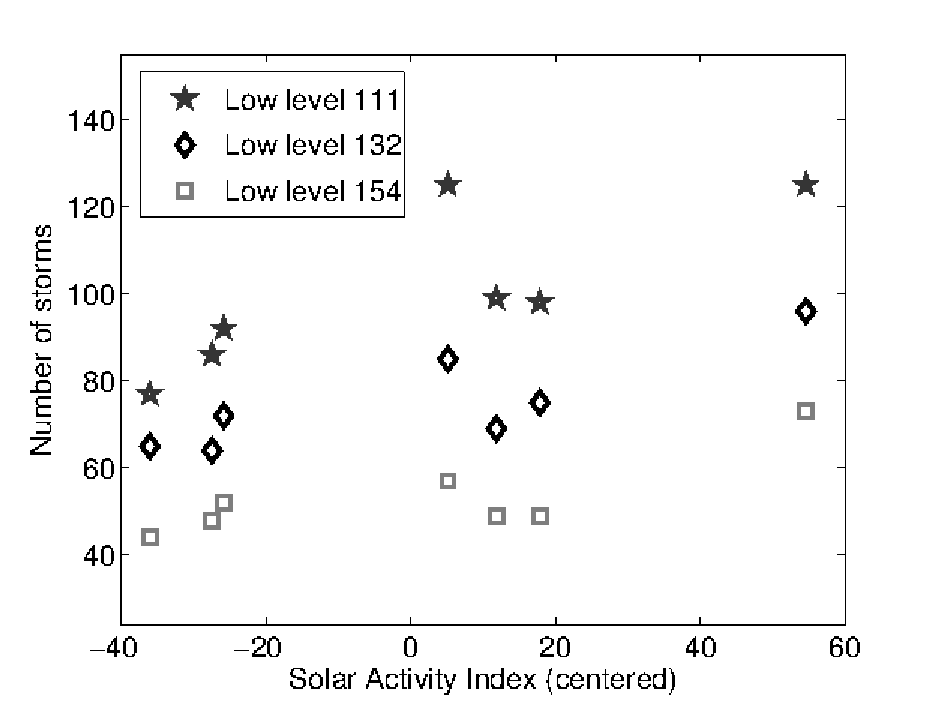}
\caption{Total number of storms per cycle for a low level of 111, 132 or 154 as a function of solar activity (centered).} 
\label{fig : NbO_AS_111} 
\end{figure}

\noindent With a low level of 111, estimation of $\beta$ gives: 
$\hat\beta =  0.0060\ \textrm{ with the 95\% C.I. } \ [ 0.0039  ; 0.0083].$
Although this value seems to be small, significance of $ \hat\beta$ has been demonstrated by a likelihood ratio test. The test of $\beta=0$ against $\beta=\hat\beta$ returns a p-value of $7.02 \times 10^{-7}$. Thus, solar activity index $X$ affects the number of storms occurring during a cycle.
Graphically, the influence of solar activity index on the number of storms per cycle is observable in Figure \ref{fig : NbO_AS_111}.\\

\subsection{$\widehat\lambda(t)$, $\widehat{\lambda_{400}}(t)$ and relative risk}
\label{section : extrapol}

\begin{figure}[H]
\center 
\includegraphics[width=8cm,height=6.2cm]{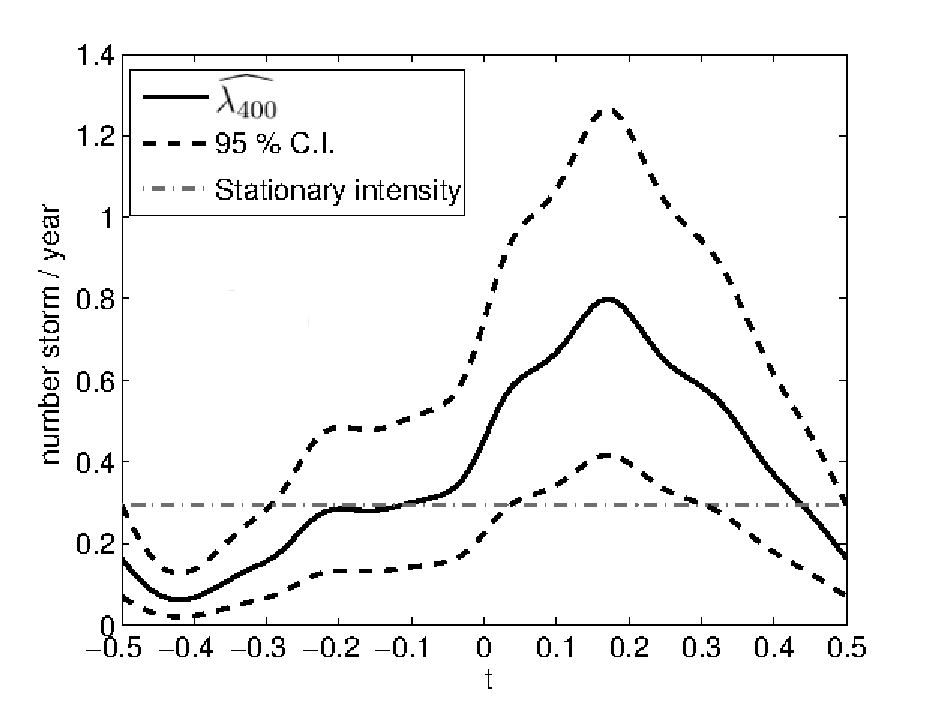}
\caption{Estimated instantaneous intensity (years$^{-1}$), with confidence interval, of the storms of level 400 obtained by extrapolation from the low level 111, for a mean solar activity of 146.7. The horizontal dash-dotted line represents the empirical frequency of storms of level 400 (stationary intensity).} 
\label{fig: intensite_111} 
\end{figure}

The extrapolation for intensity of extreme level storms is made by multiplying $\widehat\lambda(t)$ by $\widehat{P_{400}}$ (with confidence interval). The final intensity obtained, $\widehat{\lambda_{400}}$, is shown in Figure \ref{fig: intensite_111} and corresponds to the intensity of occurrence of extreme storms for a solar cycle with a mean solar activity of 146.7. 

Because of the properties of proportional hazard models, the occurrence intensity for a cycle of solar activity 146.7 corresponds to a risk factor of 1. Therefore, for a given solar cycle $j_0$ with a risk factor $ \exp ( \beta X_{j_0}) $, the relative intensity with respect to a solar cycle of solar activity 146.7 is:
$$
\lambda_{400,{j_0}}(t) = \lambda_0(t) \times P_{400} \times \exp(\beta(X_{j_0}-146.7))
$$

Using $\hat\beta$ and $\widehat{\lambda_0}$, the risk factor of a cycle can be computed and the relative intensity evaluated. For example, compared to the average level of solar activity (146.7), a cycle with a high solar activity of 180 has a risk factor of $\exp(\hat\beta(X_{j_0}-146.7))=\exp(0.0060 \times (180-146.7)) = 1.22$ with the $95\%$ C.I. $[1.14 ; 1.32]$. This implies that the extreme storms occurrence probability is expected to be 1.22 times larger during this cycle.

\subsection{Method sensitivity}

\begin{figure}[H]
\center
\includegraphics[width=7.6cm,height=6.2cm]{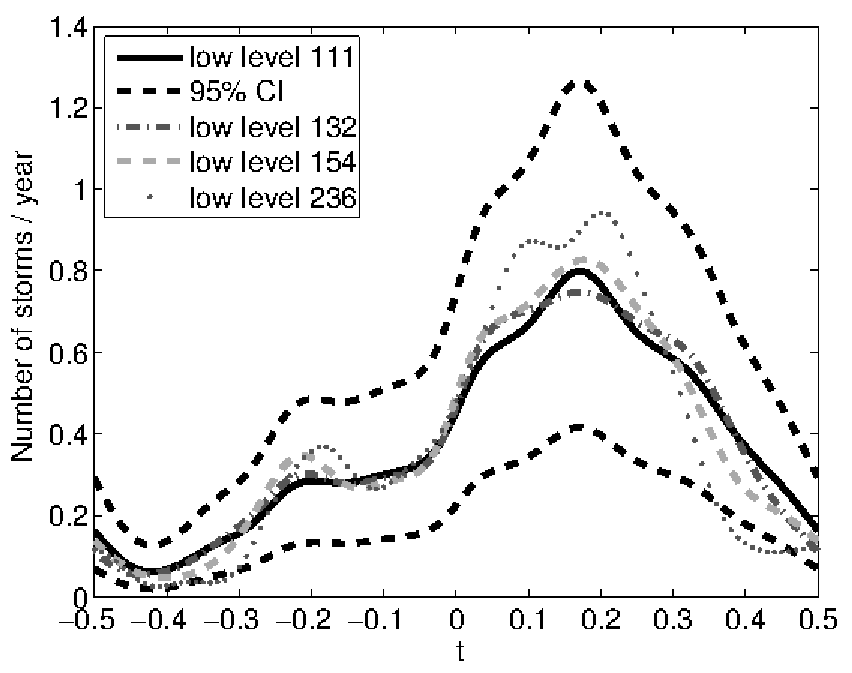}
\caption{Estimated instantaneous intensity (years$^{-1}$) of the storms of level 400 obtained by extrapolation from low levels 111 (plain line), 132 (dash-dotted line), 154 (dashed line) and 236 (dotted line).}
\label{fig: intcompare}
\end{figure}

Results presented in previous sections are given for a fixed low level (of 111). In the following, the sensitivity of the employed method is evaluated by studying its robustness to low level changes.To do so, results for three other low levels, 132, 154 and 236, are compared to that obtained with a low level of 111 (Figure \ref{fig: intcompare}). 


Figure \ref{fig: intcompare} shows that no significant difference is observable between the three curves, showing that our model is rather stable. In particular, estimated intensities for low levels 132 and 154 are always in the $95\%$ confidence interval of the intensity estimated with a low level of 111. Even if the low level is high (236) the curve shape remains coherent and most of the time within the confidence interval.

\subsection{Model extensions}
\subsubsection{Estimations on half-cycles}
Since a difference of intensity is clearly observable between the two half-cycles, a modified model is tested, where estimations are performed on each half. The number of storms for a cycle is supposed to be a non-homogeneous Poisson process, with a different intensity for each half and estimations are made separately.
However, inconclusive results are returned. Indeed, because of the presence of a normalization constant different for each half,  the occurrence intensity during the first half is higher than during the second one. Hence, this approach is not further explored.

\subsubsection{Gradient utilization}
\begin{figure}[H]
\center
\includegraphics[width=8cm,height=6.5cm]{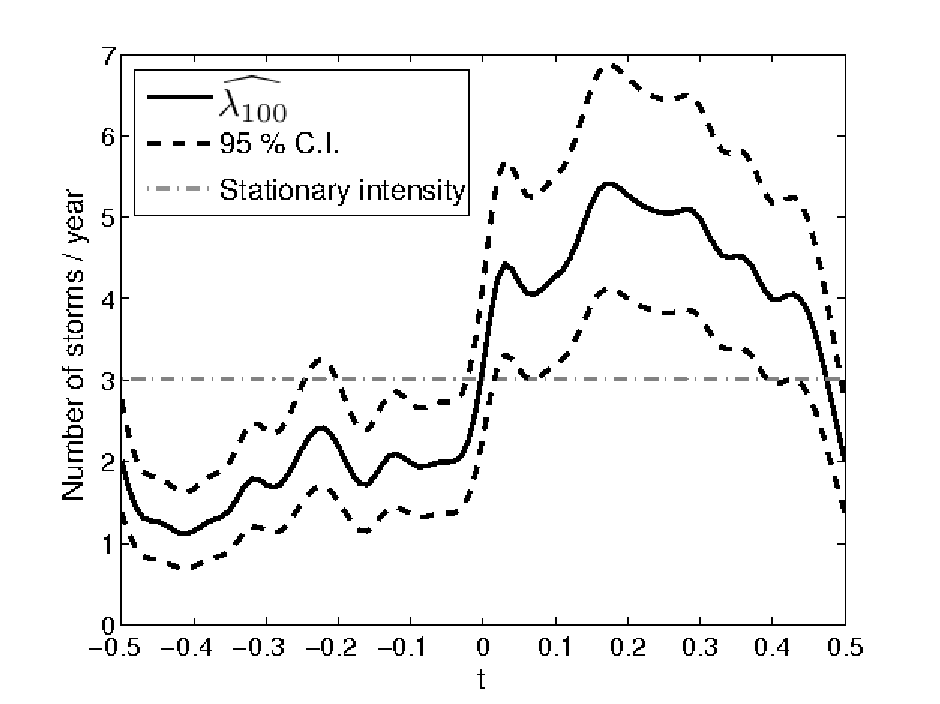}
\caption{Estimated instantaneous intensity (years$^{-1}$), with confidence interval, of the storms with extreme gradient ($\geq$ 100) obtained by extrapolation from the low gradient level 35, for a mean solar activity of 146.7. The dash-dotted line represents the empirical frequency of storms with an extreme gradient.}
\label{fig: intensite_grad}
\end{figure}

In an alternative approach, the gradient is used to characterize the strength of a storm (instead of the $ap$ index level). Gradients are computed over one time step (3 hours) and the storm gradient is defined as the maximal gradient attained during a storm. 
This approach is tested due to the occurrence of unexpected strong effects of low level storms caused by fast variations of the $ap$ index.
The same study is performed with this new definition of storm strength. The extreme gradient levels are those greater than 100 and the low one is 35. The estimation of $\beta$ gives:
$$
\hat\beta =  0.0053\ \textrm{ with the 95\% C.I. } \ [ 0.0038  ; 0.0069].
$$
Resulting values are similar to those obtained with the $ap$ index.
The estimated intensity for the storms of extreme gradient is plotted in Figure \ref{fig: intensite_grad}. The step between the two halves of the cycle is more apparent than for the previous estimation (Figure \ref{fig: intensite_111}).

Nevertheless, the gradient utilization exhibits a disadvantage. Since the $ap$ index represents a maximum over a three-hour period, the two values of $ap$ index used for the gradient computation can be separated by nearly 6 hours or only by few minutes. The real dates of these values are not known and the gradient is computed using a three-hour time step. However, the computed gradient gives an approximation of the variation speed of the $ap$ index. Moreover, since gradient is used analogously to the $ap$ index, the original model is still appropriate here.

\subsubsection{Application to the $aa$ index}

In a further sensitivity study, the proportional hazard model is applied to the $aa$ index. \\ \\
\textbf{$aa$ index description}

As the $ap$ index, $aa$ is derived from the $K$ index but for $aa$, only 2 observatories are used instead of 13  as in case of  $ap$ (one observatory is located in England and the other one in Australia). Thus, these two indices are correlated (\cite{Rifa}) but utilization of $aa$ is motivated by the larger size of the data set. Indeed, this index has been available since 1868, and 12 complete solar cycles, from the 12th (which starts on December 1878) to the 23th, are covered (against 7 for the $ap$ index). 
Especially due to its derivation process, the utilization of the $aa$ index implies adjustments; the derivation process of $aa$ is as follow:

- Firstly, for a given observatory, each magnetogram measure corresponds to a $K$ value. For example, during an extreme storm with $K = 9$, the Canberra magnetometer measurement is greater than 450 nT (nanoTesla) and the Hartland one is greater than 500 nT. Then, for each value of $K$, the scaling factor $r_K$ is found using Table 2 in \cite{calculaa} and does not depend on the observatory. The scaling factor for $K=9$ is $r_K = 667$ nT.

- Secondly, weighted average of the two rescaled values is computed. Weights depend on observatory latitude \cite{calculaa} and since locations of observatories are changing during the considered period, the weights vary too, as indicated in Table \ref{tab : coef_aa}.

\begin{table}[H]
\center
\begin{tabular}{|l|c|c|}
\hline
Period & British obs.  & Australian obs.\\
  &($w_{Eng}$ = corr. wt.)&($w_{Aus}$ = corr. wt.)\\\hline
  \hline
1868-1919   & Greenwich  & Melbourne \\
            &  (1.007)   & (0.967)   \\\hline
1920-1925   & Greenwich  & Toolangi  \\
            &  (1.007)   & (1.033)   \\\hline
1926-1956   & Albinger   & Toolangi  \\
            &  (0.934)   &  (1.033)  \\\hline
1957-1979   & Hartland   & Toolangi  \\
            &  (0.934)   &  (1.033)  \\\hline
1980-present & Hartland  & Canberra  \\
             &  (1.059)  &  (1.084)  \\\hline
\end{tabular}
\caption{Sources of $aa$ index for different time intervals and corresponding weights.}
\label{tab : coef_aa}
\end{table}

- Finally, the mathematical expression of $aa$ is:
\begin{eqnarray}
aa = 1/2 [ w_{Eng} . r_{K_{Eng}} + w_{Aus} . r_{K_{Aus}}].
\label{eq: aa}
\end{eqnarray}
For an extreme storm with $K=9$ measured after 1980 the corresponding $aa$ index is $aa = 1/2 [1.059\times 667 + 1.084 \times 667] \simeq 714.69 \textrm{ rounded to } 715$ (the maximal value of $aa$).\\ \\
\textbf{Definition of an extreme storm}

The previous relation raises the question of rescaling low and extreme levels used for estimation. A particular attention is taken to the change of magnetic observatories over time, in both England and Australia. These changes are handled by creating two different levels for the periods 1868-1956 (cycles 12 to 18) and 1957-present (cycles 19 to 23) since for each period, weights are quite equivalent. In order to detect all the extreme storms, the lower level has been selected for each period. 

An extreme storm is characterised by an $ap$ index of 400 or a $K$ index of 9 according to Table \ref{tab : relationKap}. When the corresponding extreme level for $aa$ is computed using Equation (\ref{eq: aa}), only 3 storms are detected as extreme since 1878 (the beginning of the first complete solar cycle) and only 2 since 1932 (the period when the $ap$ index is also available). 23 extreme storms are counted with the $ap$ index. These numbers are too different and comparison is not achievable.

Since the $aa$ index is less global than the $ap$ index, the following rule is adopted: considering the time period  when the two indices are available, an extreme storm for the $ap$ index should also be detected as extreme for the $aa$ index, even if additional local events are detected. Under this condition, the $aa$ extreme level of 500 seems to be adapted with 23 extreme storms detected for $aa$ against 23 for the $ap$ index (only the second period level is given, the first is automatically computed and is 472). Thus, an extreme storm for $aa$ is defined as a storm of level greater or equal to 500 and the total number of extreme storms detected for the $aa$ index is 33.\\ \\
\textbf{Results}

The model is then applied to the $aa$ data set with an extreme level of 500, a low level of 130 and a run length of 2. The low level selection is less important than the extreme one because of the model non-sensitivity to this parameter, that remains valid in this application. The low level of 130 is selected similarly to the extreme level.
The solar activity index covariate $X$ is now $X = $(74.6,\ 87.9,\  64.2,\  105.4,\  78.1,\  119.2,\  151.8,\  201.3,\  110.6,\  164.5,\  158.5,\  120.8) with a mean solar activity of about 119.7.
The estimations for $\beta$ and $P_{500}$ are:
$$
\hat\beta =  0.0068\ \textrm{ with the 95\% C.I. } \ [ 0.0055  ;   0.0081],
$$
$$
\widehat{P_{500}} =  2.49 \% \ \textrm{ with the 95\% C.I. } \ [ 1.65 ; 3.33].
$$

\begin{figure}[H]
\center
\includegraphics[width=8cm,height=6.3cm]{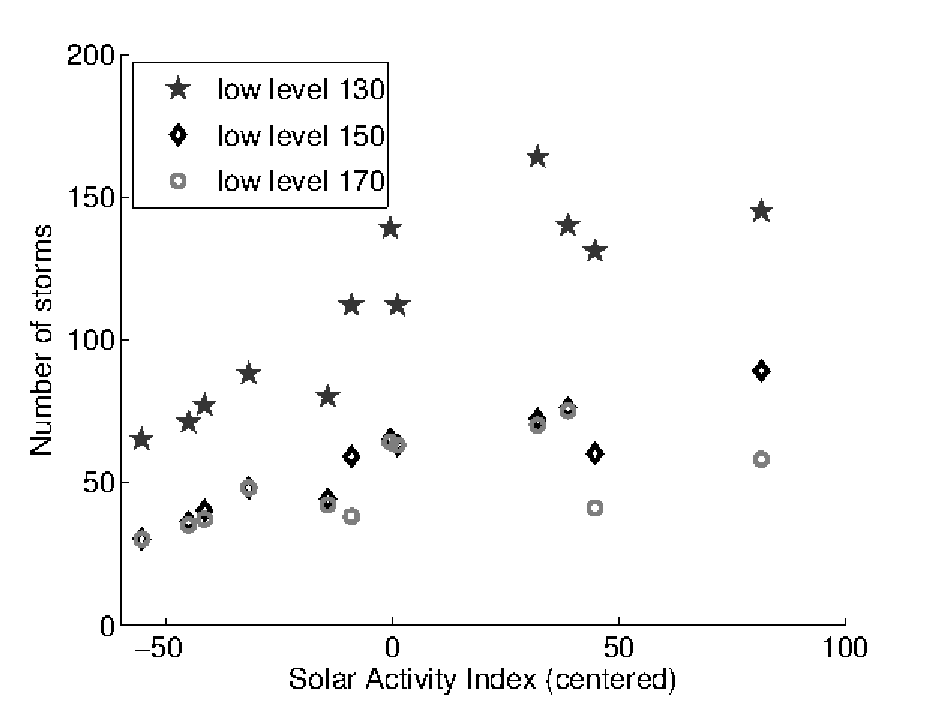}
\caption{Total number of storms per cycle for low levels of 130, 150 and 170 according to the solar activity (centered).} 
\label{fig : NbO_AS_aa} 
\end{figure}

The $\hat\beta$ value is in the $95 \%$ confidence interval of the estimation using $ap$. Since more cycles are available, the confidence interval has been now shortened. Influence of solar activity is a little bit more apparent (Figure \ref{fig : NbO_AS_aa}).

\begin{figure}[H]
\center
\includegraphics[width=7.9cm,height=6cm]{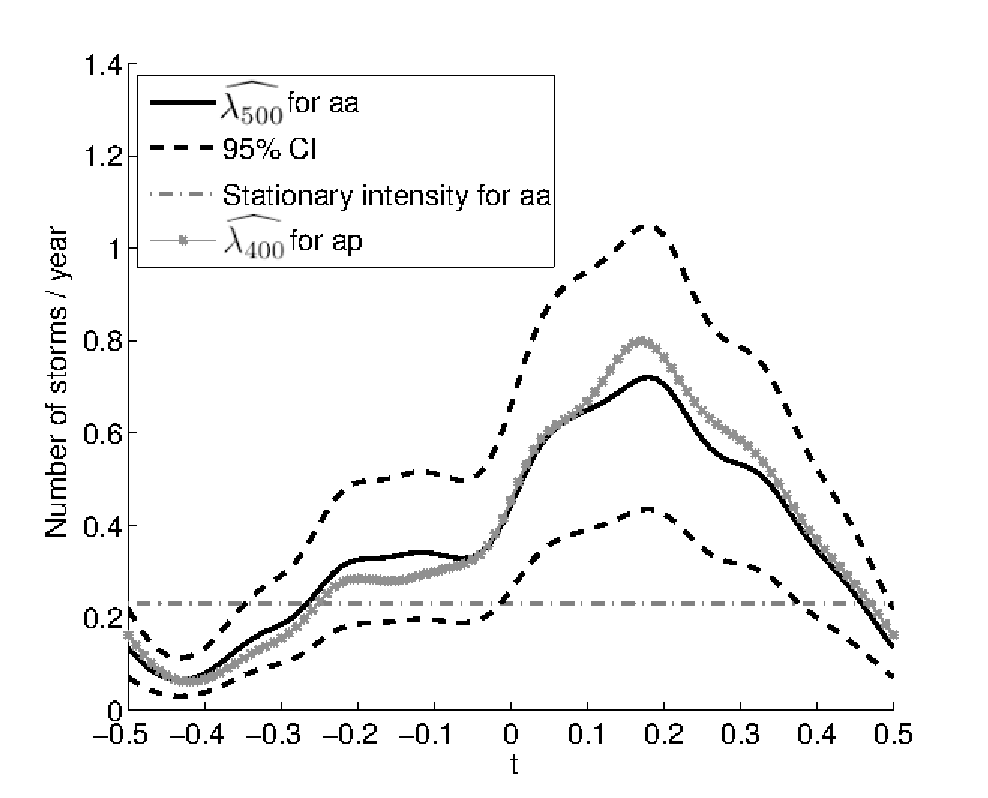}
\caption{Estimated instantaneous intensity (years$^{-1}$) of the storms defined by $aa \geq 500$ (plain line). Estimation is obtained by extrapolation from the low level 130, for a mean solar activity of 119.7. The dotted lines represent the $95 \%$ confidence interval. The gray horizontal line corresponds to the empirical frequency of extreme storms. The gray curve with stars is the intensity estimated with the $ap$ data set (for a low level of 111).} 
\label{fig : int_aa_ap2} 
\end{figure}

The final intensity estimation $\widehat{\lambda_{500}}$ is plotted in Figure \ref{fig : int_aa_ap2}.

The empirical frequency of extreme storms is now 0.23 (since 1878), against 0.29 for the $ap$ data set. This decrease can be explained by addition of solar cycles with a lower solar activity index (the mean solar activity since 1878 is 119.7, against 146.7 for the period 1932 to present). 

\begin{figure}[H]
\center
\includegraphics[width=7.9cm,height=5.5cm]{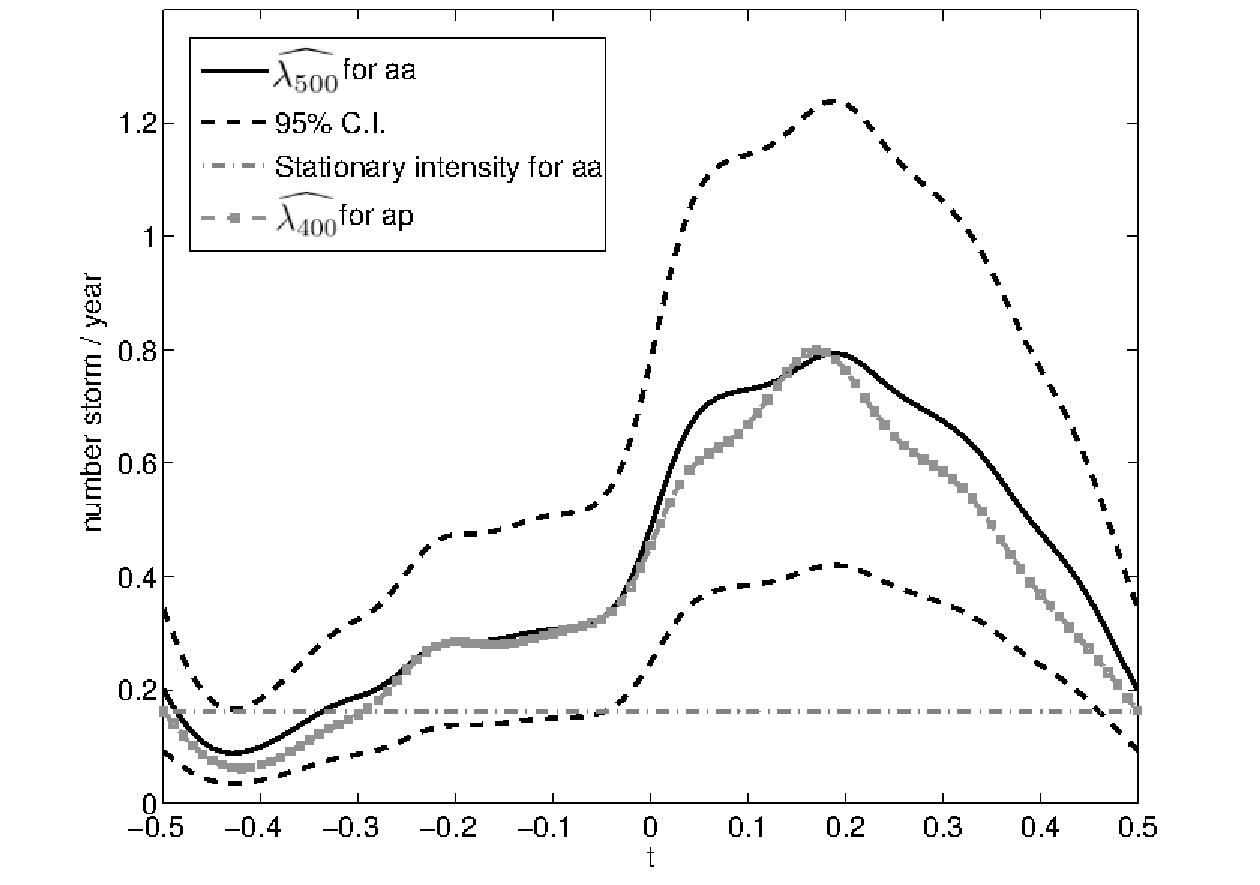}
\caption{Similar to Figure \ref{fig : int_aa_ap2} with the $aa$ index considered on the 7 complete cycles of the $ap$ data set.} 
\label{fig : int_aa32} 
\end{figure}

Figure \ref{fig : int_aa_ap2} should be read carefully since intensity curves correspond respectively to a theoretical cycle with a mean solar activity of 146.7 for the $ap$ index and 119.7 for the $aa$ index. Results using the $aa$ index only for the $ap$ index period of availability are also computed (and then correspond to the same mean solar activity) and plotted in Figure \ref{fig : int_aa32}.
The most important result of this paragraph is the similarity of the two intensity curves in Figure \ref{fig : int_aa32}, showing that the method is stable against a change of data set.

\subsection{Prevision}

The model proposed here also can be used for forecasting the intensity of extreme storms for the current 24th solar cycle. Data needed are:\\
- beginning date;\\
- solar activity index (the $X$ covariate);\\
- peak date;\\
- end date, and so, cycle length $D$.

\begin{figure}[H]
\center
\includegraphics[width=8.5cm,height=7cm]{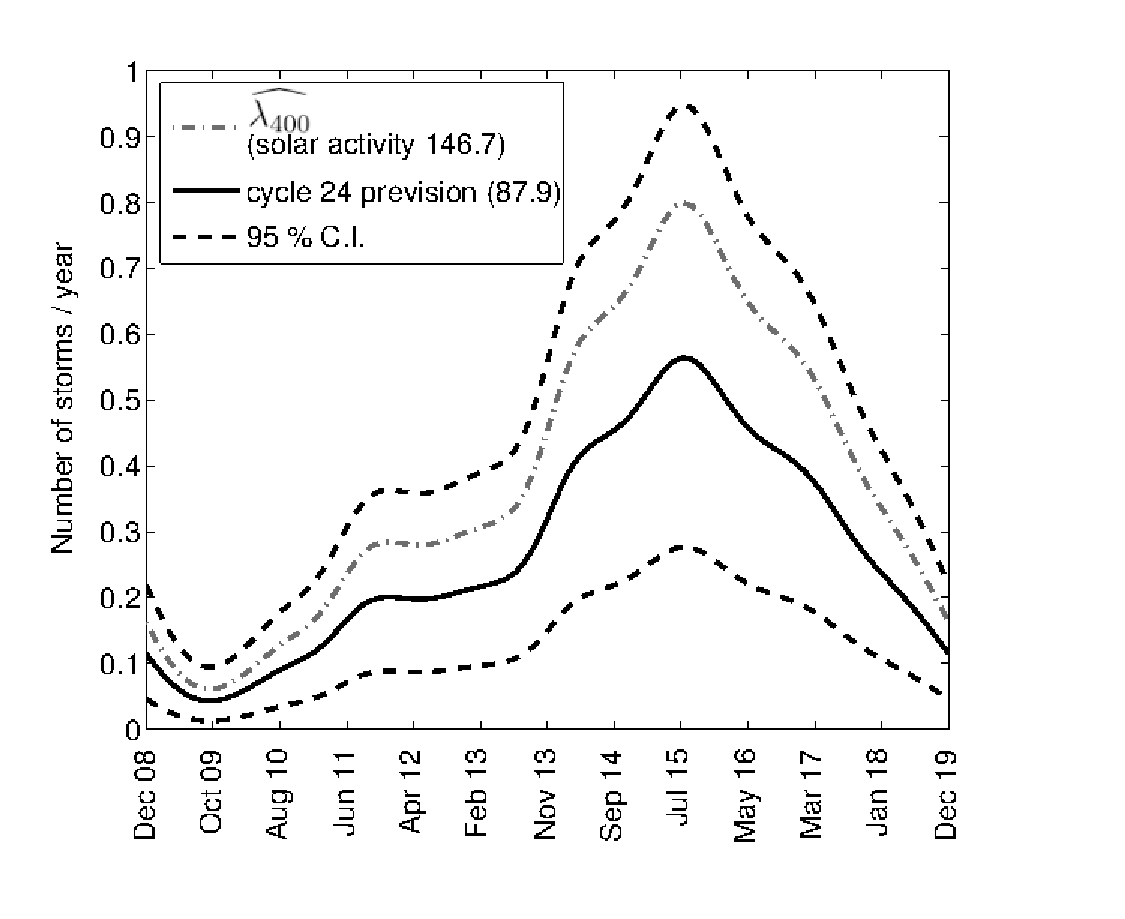}
\caption{Estimation and prediction of instantaneous intensity (years$^{-1}$) of the storms of level 400 for the 24th solar cycle, with confidence interval. For comparison, in dash dotted gray, same intensity for a cycle with a mean solar activity index of 146.7.}
\label{fig: prev24}
\end{figure}

Beginning date is known but the three other data must be estimated since the 24th cycle has not ended. NOAA predictions are used and presented in Table \ref{tab: covar}.
Estimation (from the beginning to present (April 2013)) and prediction are represented in Figure \ref{fig: prev24} (plain line). Solar activity index of the current cycle is estimated at 87.9 and the mean solar activity index of the  previous seven cycles is 146.7. Thus the risk factor of the 24th cycle (compared to a cycle with mean solar activity of 146.7) is $\exp(\hat\beta(87.9-146.7)) = 0.70$. Thus, occurrence intensity for the current cycle is  0.70 times weaker than the one estimated for a cycle of mean solar activity (gray dash-dotted line in Figure \ref{fig: prev24}).

\section{CONCLUSION}

In this paper, an ionospheric storm occurrence risk estimation is performed. As extreme storms are rare events, classical methods for probability estimation can not be applied. Moreover, the occurrence probability of ionospheric storms varies over time and the estimation method should take into account this non-stationarity. Based on the Cox model, the construction of an innovative approach using a proportional hazard model has been motivated by the proportional hazard assumptions. In this model, the intensity function can depend on time and a covariate influence can be handled. Mathematically, the number of storms for a cycle is modeled by a non-homogeneous Poisson process of whom the intensity depends on time and one covariate (the solar activity index).

However, extreme events are too scarce to apply straightforward estimations. Then, adopting an Extreme Value Theory process, estimations are made using high events, less rare, and the occurrence intensity across a solar cycle is extrapolated to extreme events, as in the classical Peaks Over Threshold method.

The $ap$ index is chosen to describe ionospheric perturbations. Advantages of this index are its global character, the large amount of data (7 complete solar cycles) and the completeness of time series. 

Sensitivity of the model is investigated in two different ways:

- firstly, the model response to a change of low level parameter is examined (this level is used to define the high level storms). Analysis are performed for different low levels. First results of $ap$ index are given for a low level of 111 and a comparison is made using three other low levels: 132, 154 and 236. The shape similarity between all the intensity curves confirms the stability of the method.

- Secondly, the data set used for estimation is changed. In the first approach, the maximal $ap$ value reached during a storm is used to characterise its strength. A second approach considers gradients of storms to classify them as high or extreme events. In a third step, an analogous study is performed using the $aa$ index. All the results obtained are coherent, showing a strong similarity between all the estimated intensity curves.

The main results of this study show a strong correlation between the occurrence intensity of magnetic storms and their positions on the solar cycle. This intensity is higher during the second half of the cycle. 

This work also shows that solar activity has an influence on occurrence intensity of extreme storms.
With the proportional hazard model, the occurrence intensity of a cycle can be expressed as a relative risk compared to a cycle with a mean solar activity index (146.7 for the $ap$ data set and 119.7 for $aa$). The proposed model can also be used to forecast the intensity of occurrence of extreme storms for the current 24th solar cycle.

\appendix
\section*{Appendix}
\section{Chi-square test}
\label{annexe: chi2}

The utilization of the parameter $\widehat{P_{400}}$ (the probability for a high level storm to grow into an extreme level storm) for the extrapolation supposes that the level reached by a storm does not depend on its   occurrence time. $\widehat{P_{400}}$ is then estimated using the empirical frequency. To verify this independence assumption, a Chi-square test is performed.

The warped time interval $[-0.5,0.5]$ is partitioned into two parts of low and high intensity using the final estimated instantaneous intensity $\widehat{\lambda_{400}}$. The period of low intensity is defined as $T_l = \{ t \ | \widehat{\lambda_{400}}(t) < 0.30\}$ and the period of high intensity as $T_h = \{ t \ | \widehat{\lambda_{400}}(t) \geq 0.30\}$. A threshold of 0.30 is the lowest possible which makes the application of a Chi-square test possible (under the hypothesis of independence the expected number of storms in each cell should be greater or equal to 5). 

The numbers of high and extreme storms for $T_h$ and $T_l$ are gathered in a $2 \times 2$ contingency table (Table \ref{tab : conting}). 

\begin{table}[H]
\center
\begin{tabular}{|c|c|c|}
\cline{2-3} \multicolumn{1}{c|}{} & Extreme storms  & High storms\\ \hline
$T_h$ &  19 & 538 \\ \hline
$T_l$ & 4 & 164  \\ \hline
\end{tabular}\\
\caption{Number of extreme and high storms for low and high intensity periods, respectively.}
\label{tab : conting}
\end{table}

The Chi-square statistic is computed for this contingency table and the corresponding p-value (for a significance level of $\alpha = 5\%$) is 0.50. Thus, the independence hypothesis is not rejected.
The same test is applied with different thresholds for the warped time partition into $T_h$ and $T_l$ (0.40, 0.50 and 0.60) and always leads to the same conclusion with p-values of 0.67, 0.39 and 0.61, respectively.

\section{Maximum likelihood estimator of $\beta$}
\label{annexe}
The use of $N_j$ (the total number of high level storms in cycle $j$) instead of $N_j(t)$ for the estimation of $\beta$ raises the question of sufficiency of this statistic. Consider only one cycle and the model:
$$
N(t) \sim  \mathcal P \left(  \lambda_0 (t)dt  \ D \exp(\beta X) \right),\quad \textrm{ for } t \in [-0.5 , 0.5 ].
$$
Then, consider $\Delta_1, \Delta_2, ..., \Delta_n$ a partition of [-0.5, 0.5] into $n$ sub-segments. For $i=1...n$, note $N(\Delta_i)= \int_{\Delta_i} dN(t)$ the number of events in $\Delta_i$. The $\{N(\Delta_i), i=1...n\}$ are independent variables (by definition of a Poisson process) and $N(\Delta_i) \sim \mathcal P \left(  \left[ \int_{\Delta_i} \lambda_0 (s)ds \right]  \ D \ \exp(\beta X) \right)$. Noting $C_i = \int_{\Delta_i} \lambda_0 (s)ds  \ D$, the Log-likelihood with respect to the counting measure (in which we integrate the weights $1/N(\Delta_i)!\ $) is:
$$
- \exp(\beta X) \sum_{i=1}^n  C_i + \sum_{i=1}^n [ N(\Delta_i) \log(C_i)] + \beta X \sum_{i=1}^n N(\Delta_i).
$$
$\beta$ is linked to the $N(\Delta_i)$ only by the term $\sum_{i=1}^n N(\Delta_i)$. Hence there is no loss of information when using the total number of events per cycle for the estimation of $\beta$.

It is now possible to compute the maximum likelihood estimator. For the $j^{th}$ cycle, the likelihood with respect to the counting measure with weights $1/N_j! \, $ is, noting $\alpha = \int_{-1/2}^{1/2} \lambda_0 (s)ds$
$$
\exp \left(   -\alpha \, D_j \exp(\beta X_j)  \right) (\alpha \, D_j \exp(\beta X_j))^{N_j},
$$
and the Log-likelihood for all the $J$ cycles:
\begin{eqnarray*}
-\alpha \sum_{j=1}^J  D_j \exp(\beta X_j) \! &+&\log(\alpha)  \sum_{j=1}^J N_j \\
&+& \! \sum_{j=1}^J N_j \log(D_j) + \beta \sum_{j=1}^J N_jX_j.
\end{eqnarray*}  

The derivatives in $\alpha$ anb $\beta$ respectively give:
$$
\sum_{j=1}^J  D_j \exp(\beta X_j) = \frac{\sum_{j=1}^J N_j}{\alpha},
$$
and
$$
\alpha \sum_{j=1}^J  D_j X_j  \exp(\beta X_j) = \sum_{j=1}^J N_j X_j \, .
$$
Replacing $\alpha$ by the solution of the first equation, we obtain:
$$
\sum_{j=1}^J  D_j X_j  \exp(\beta X_j) \sum_{j=1}^J N_j = \sum_{j=1}^J  D_j  \exp(\beta X_j)  \sum_{j=1}^J N_j X_j \, .
$$
This implicit equation can be resolved only numerically (by the secant method).

The Fisher information matrix is also computable:
$$
\left( \begin{array}{cc}
\alpha^{-1} \sum_{j=1}^J D_j \exp(\beta X_j) & \sum_{j=1}^J D_j X_j \exp(\beta X_j)  \\
\sum_{j=1}^J D_j X_j \exp(\beta X_j) & \alpha \sum_{j=1}^J D_j X^2_j \exp(\beta X_j)  \\
\end{array} \right).
$$
The (2,2) coefficient of the inverse matrix of the Fisher information matrix provides the variance of $\hat \beta$, used for the construction of a confidence interval.

\section*{Acknowledgements}
The authors wish to thank the editor, the associate editor and two referees of Statistics and Its Interface (SII) for providing helpful comments and suggestions. This PhD candidate work at IMT lab of UPS is founded by CNES and TASF.

\bibliographystyle{plain}
\bibliography{biblio}

\begin{thebibliography}{10}

\bibitem{noaaKpap}
\url{http://www.ngdc.noaa.gov/stp/GEOMAG/kp_ap.html}.

\bibitem{nasaSSN}
\url{http://www.sidc.be/sunspot-data/}.

\bibitem{report2008}
Severe space weather events--understanding societal and economic impacts.
\newblock {\em Workshop report of the Societal and Economic Impacts of Severe
  Space Weather Events}, 2008.

\bibitem{noaaPrev24}
\url{http://www.swpc.noaa.gov/ftpdir/weekly/Predict.txt}, Update on April 8,
  2013, consulted on April 8, 2013.

\bibitem{Aalen2008}
O.O. Aalen, Ø. Borgan, and H.F. Gjessing.
\newblock {\em Survival and event history analysis. A process point of view}.
\newblock Springer, New York, 2008.

\bibitem{AndersenGill}
P.K. Andersen and R.D. Gill.
\newblock Cox's regression model for counting processes: a large sample study.
\newblock {\em Ann. Statist.}, 10(4):1100--1120, 1982.

\bibitem{Anderson1970}
C.~W. Anderson.
\newblock Extreme value theory for a class of discrete distributions with
  applications to some stochastic processes.
\newblock {\em J. Appl. Probability}, 7:99--113, 1970.

\bibitem{Bowman1984}
A.W. Bowman.
\newblock An alternative method of cross-validation for the smoothing of
  density estimates.
\newblock {\em Biometrika}, 71(2):353--360, 1984.

\bibitem{Coles2001}
S.~Coles.
\newblock {\em An introduction to statistical modeling of extreme values}.
\newblock Springer Series in Statistics. Springer-Verlag London Ltd., London,
  2001.

\bibitem{CooleyNaveau2007}
D.~Cooley, D.~Nychka, and P.~Naveau.
\newblock Bayesian spatial modeling of extreme precipitation return levels.
\newblock {\em J. Am. Statist. Ass.}, 102(479):824--840, 2007.

\bibitem{Cox72}
D.R. Cox.
\newblock Regression models and life-tables.
\newblock {\em J. R. Statist. Soc. B}, 34:187--220, 1972.

\bibitem{Hall1983}
P.~Hall.
\newblock Large sample optimality of least squares cross-validation in density
  estimation.
\newblock {\em Ann. Statist.}, 11(4):1156--1174, 1983.

\bibitem{Husler86}
J.~Hüsler.
\newblock Extreme values of non-stationary random sequences.
\newblock {\em J. Appl. Probability}, 23(4):937--950, 1986.

\bibitem{Jonathan2011}
P.~Jonathan and K.~Ewans.
\newblock Modeling the seasonality of extreme waves in the gulf of mexico.
\newblock {\em J. Offshore Mech. Arct. Eng.}, 133(2):113--123, 2011.

\bibitem{calculaa}
J.J. Love.
\newblock Long-term biases in geomagnetic k and aa indices.
\newblock {\em Ann. Geophys.}, 29:1365--1375, 2011.

\bibitem{NielsenLinton}
J.P. Nielsen and O.B. Linton.
\newblock Kernel estimation in a nonparametric marker dependent hazard model.
\newblock {\em Ann. Statist.}, 23(5):1735--1748, 1995.

\bibitem{Rifa}
E.~Rifa.
\newblock Études des relations entre indices solaires et géomagnetiques.
\newblock {\em Intership report}.

\bibitem{Schwabe}
H.~Schwabe.
\newblock Solar observations during 1843.
\newblock {\em Astronomische Nachrichten}, 20(495):233--236, 1843.

\bibitem{Norbert}
N.~Suard, J.M. Azaïs, S.~Gadat, D.~Debailleux, and E.~Rifa.
\newblock Assessment of an ionosphere storm occurrence risk.
\newblock {\em European Navigation Conference}, 2011.

\bibitem{NOAAmemo}
M.~Weaver, W.~Murtagh, C.~Balch, D.~Biesecker, and L.~Combs.
\newblock Halloween space weather storms of 2003.
\newblock {\em NOAA Technical Memorandum OAR SEC-88}, June 2004.

\end{thebibliography}

\end{document}